\documentclass[12pt]{amsart}
\usepackage{graphicx}
\usepackage{amssymb}
\usepackage{amsfonts}
\usepackage{amsmath}
\usepackage{array}
\usepackage{rotating}

\usepackage[matrix,frame]{xypic}
\newdimen\vcadre\vcadre=0.1cm % marges verticales de la boite
\newdimen\hcadre\hcadre=0.1cm % marges horizontales de la boite
\def\GrTeXBox#1{\vbox{\vskip\vcadre\hbox{\hskip\hcadre%
      $#1$%
   \hskip\hcadre}\vskip\vcadre}}
\def\arx#1[#2]{\ifcase#1 \relax \or%
  \ar @{-}[#2]  \or%
  \ar @2{-}[#2] \or%
  \ar @{--}[#2] \or%
  \ar @2{.}[#2] \or%
  \ar @{~}[#2]  \fi}

\usepackage{tikz}
\usetikzlibrary{arrows.meta}

\newtheorem{example}{Example}[section]
\newtheorem{note}[example]{Note}
\newtheorem{remark}[example]{Remark}
\newtheorem{theorem}[example]{Theorem}

\newtheorem{corollary}[example]{Corollary}

\newtheorem{definition}[example]{Definition}
\newtheorem{proposition}[example]{Proposition}

\newtheorem{lemma}[example]{Lemma}

\def\Proof{\noindent \it Proof -- \rm}
\def\qed{\hspace{3.5mm} \hfill \vbox{\hrule height 3pt depth 2 pt width 2mm}
\bigskip}

\def\XX{{\bf X}}

\def\FQSym{{\bf FQSym}}
\def\WQSym{{\bf WQSym}}
\def\WSym{{\bf WSym}}

\def\CQSym{{\bf CQSym}}

\def\ev{{\rm ev}}

\def\std{{\rm std}}

\def\st{{\rm st}}

\def\<{\langle}
\def\>{\rangle}

\def\N{{\bf N}}
\def\NN{{\mathbb N}}
\def\K{{\bf K}}
\def\KK{{\mathbb K}}

\def\pack{\operatorname{pack}}

\def\F{{\bf F}}

\def\G{{\bf G}}
\def\M{{\bf M}}

\def\SG{{\mathfrak S}}

\def\Sym{{\bf Sym}}
\def\QSym{{\it QSym}}

\def\Des{\operatorname{Des}}

\def\maj{{\rm maj}}

\def\PW{{\rm PW}}

\def\GPG{{\mathcal{G}}}
\def\GPD{{\mathcal{D}}}

\def\X{{\rm X}}
\def\LLT{{\rm LLT}}
\def\asc{{\rm asc}}
\def\shuff#1#2{\mathbin{
\hbox{\vbox{ \hbox{\vrule \hskip#2 \vrule height#1 width 0pt
}%
\hrule}%
\vbox{ \hbox{\vrule \hskip#2 \vrule height#1 width 0pt
\vrule }%
\hrule}%
}}}

\def\shuf{{\mathchoice{\shuff{7pt}{3.5pt}}%
{\shuff{6pt}{3pt}}%
{\shuff{4pt}{2pt}}%
{\shuff{3pt}{1.5pt}}}}%
\def\shuffle{\,\shuf\,}
\def\inv{{\rm inv}}

\def\ogf#1#2{\genfrac{}{}{0pt}{}{#1}{#2}}

\makeatletter
\newsavebox{\@brx}
\newcommand{\llangle}[1][]{\savebox{\@brx}{\(\m@th{#1\langle}\)}%
  \mathopen{\copy\@brx\kern-0.5\wd\@brx\usebox{\@brx}}}
\newcommand{\rrangle}[1][]{\savebox{\@brx}{\(\m@th{#1\rangle}\)}%
  \mathclose{\copy\@brx\kern-0.5\wd\@brx\usebox{\@brx}}}
\makeatother

\def\Tabvrule{\vrule width-0.4pt}       % Difference de largeur
\def\Tabhrule{\hrule \hrule height-0.4pt} % Difference de hauteur
\def\Tabstrut{\vrule height2.2ex % Sur la ligne
                     depth0.8ex  % Sous la ligne
                     width0ex    % centrage horizontal
\relax}

\def\PasCase#1{\omit%
            $\vcenter{\hbox {\vbox to 0.4pt{}}
               \hbox{\makebox[3ex]{\Tabstrut$#1$}}}%
               \Tabvrule$}
\def\PasCasePoint{\PasCase{\cdot}}
\def\DessinCarre#1{%
    \vcenter{\hbox{}\hrule
             \hbox{\vrule\makebox[3ex]{\Tabstrut$#1$}\vrule}\Tabhrule}%
             \Tabvrule}
\def\GenRuban#1{\vcenter{\halign{&$\DessinCarre{##}$\cr#1}}\egroup}

\def\sTabvrule{\vrule width-0.4pt}
\def\sTabhrule{\hrule \hrule height-0.4pt}
\def\sTabstrut{\vrule height1.6ex depth0.6ex width0ex \relax}
\def\sDessinCarre#1{%
    \vcenter{\hbox{}\hrule
             \hbox{\vrule\makebox[2.3ex]%
                  {\sTabstrut$\scriptstyle#1$}\vrule}\sTabhrule}%
             \sTabvrule}
\def\sGenRuban#1{\vcenter{\halign{&$\sDessinCarre{##}$\cr#1}}\egroup}

\def\ruban{%
  \bgroup
  \let\ =\omit
  \let\\=\cr
  \let\x=\times
  \let\.=\PasCasePoint
  \offinterlineskip
  \GenRuban}

\def\sruban{%
  \bgroup
  \let\ =\omit
  \let\x=\times
  \let\\=\cr
  \offinterlineskip
  \sGenRuban}
\newdimen\Squaresize \Squaresize=14pt
\newdimen\Thickness \Thickness=0.5pt

\def\Square#1{\hbox{\vrule width \Thickness
   \vbox to \Squaresize{\hrule height \Thickness\vss
      \hbox to \Squaresize{\hss#1\hss}
   \vss\hrule height\Thickness}
\unskip\vrule width \Thickness}
\kern-\Thickness}

\def\Vsquare#1{\vbox{\Square{$#1$}}\kern-\Thickness}

\def\young#1{
\vbox{\smallskip\offinterlineskip
\halign{&\Vsquare{##}\cr #1}}}

\def\boxit#1#2{\setbox1=\hbox{\kern#1{#2}\kern#1}%
\dimen1=\ht1 \advance\dimen1 by #1 \dimen2=\dp1 \advance\dimen2 by #1
\setbox1=\hbox{\vrule height\dimen1 depth\dimen2\box1\vrule}%
\setbox1=\vbox{\hrule\box1\hrule}%
\advance\dimen1 by .4pt \ht1=\dimen1
\advance\dimen2 by .4pt \dp1=\dimen2 \box1\relax}

\def\raff{\geqslant_{\rm ref}}

\def\PC{\rm PC\,}

\def\taille{.5}

\def\m{{\bf m}}
\def\mt{{\bf \tilde m}}
\def\dn{{\rm dn}}
\def\vPhi{\check{\Phi}}
\def\bLLT{{\bf LLT}}
\def\DST{{\rm DST\,}}

\def\gunun{
\begin{tikzpicture}
\begin{scope}[every node/.style={circle,scale=.5,fill=white,draw}]
    \node (A) at (0,0) {};
\end{scope}
\end{tikzpicture}
}

\def\gdeuxun{
\begin{tikzpicture}
\begin{scope}[every node/.style={circle,scale=.5,fill=white,draw}]
    \node (A) at (0,0) {};
    \node (B) at (1*\taille,0) {};
\end{scope}
\end{tikzpicture}
}

\def\gdeuxdeux{
\begin{tikzpicture}
\begin{scope}[every node/.style={circle,scale=.5,fill=white,draw}]
    \node (A) at (0,0) {};
    \node (B) at (1*\taille,0) {};
\end{scope}

\begin{scope}[>={Stealth[black]},
              every edge/.style={draw=black,thick}]
    \path [-] (A) edge (B);
\end{scope}
\end{tikzpicture}
}

\def\gtroisun{
\begin{tikzpicture}
\begin{scope}[every node/.style={circle,scale=.5,fill=white,draw}]
    \node (A) at (0,0) {};
    \node (B) at (1*\taille,0) {};
    \node (C) at (2*\taille,0) {};
\end{scope}
\end{tikzpicture}
}

\def\gtroisdeux{
\begin{tikzpicture}
\begin{scope}[every node/.style={circle,scale=.5,fill=white,draw}]
    \node (A) at (0,0) {};
    \node (B) at (1*\taille,0) {};
    \node (C) at (2*\taille,0) {};
\end{scope}

\begin{scope}[>={Stealth[black]},
              every edge/.style={draw=black,thick}]
    \path [-] (A) edge (B);
\end{scope}
\end{tikzpicture}
}

\def\gtroistrois{
\begin{tikzpicture}
\begin{scope}[every node/.style={circle,scale=.5,fill=white,draw}]
    \node (A) at (0,0) {};
    \node (B) at (1*\taille,0) {};
    \node (C) at (2*\taille,0) {};
\end{scope}

\begin{scope}[>={Stealth[black]},
              every edge/.style={draw=black,thick}]
    \path [-] (B) edge (C);
\end{scope}
\end{tikzpicture}
}

\def\gtroisquatre{
\begin{tikzpicture}
\begin{scope}[every node/.style={circle,scale=.5,fill=white,draw}]
    \node (A) at (0,0) {};
    \node (B) at (1*\taille,0) {};
    \node (C) at (2*\taille,0) {};
\end{scope}

\begin{scope}[>={Stealth[black]},
              every edge/.style={draw=black,thick}]
    \path [-] (A) edge (B);
    \path [-] (B) edge (C);
\end{scope}
\end{tikzpicture}
}

\def\gtroiscinq{
\begin{tikzpicture}
\begin{scope}[every node/.style={circle,scale=.5,fill=white,draw}]
    \node (A) at (0,0) {};
    \node (B) at (1*\taille,0) {};
    \node (C) at (2*\taille,0) {};
\end{scope}

\begin{scope}[>={Stealth[black]},
              every edge/.style={draw=black,thick}]
    \path [-] (A) edge (B);
    \path [-] (B) edge (C);
    \path [-] (A) edge[bend left=60] (C);
\end{scope}
\end{tikzpicture}
}

\def\gsixex{
\begin{tikzpicture}
\begin{scope}[every node/.style={circle,scale=.5,fill=white,draw}]
    \node (A) at (0,0) {};
    \node (B) at (1*\taille,0) {};
    \node (C) at (2*\taille,0) {};
    \node (D) at (3*\taille,0) {};
    \node (E) at (4*\taille,0) {};
    \node (F) at (5*\taille,0) {};
\end{scope}

\begin{scope}[>={Stealth[black]},
              every edge/.style={draw=black,thick}]
    \path [-] (A) edge (B);
    \path [-] (B) edge (C);
    \path [-] (C) edge (D);
    \path [-] (D) edge (E);
    \path [-] (E) edge (F);
    \path [-] (B) edge[bend left=60] (D);
    \path [-] (D) edge[bend left=60] (F);
\end{scope}
\end{tikzpicture}
}

\title{Noncommutative unicellular LLT polynomials}

\author[ J.-C.~Novelli and J.-Y.~Thibon]%
{Jean-Christophe Novelli and Jean-Yves Thibon}

\address[] {LIGM, Universit\'e
Gustave-Eiffel, CNRS, ENPC, ESIEE-Paris \\
5 Boulevard Descartes \\Champs-sur-Marne \\77454 Marne-la-Vall\'ee cedex 2 \\
FRANCE}
\email[Jean-Christophe Novelli]{novelli@univ-mlv.fr}
\email[Jean-Yves Thibon]{jyt@univ-mlv.fr} 



\date{}

\begin{document}

\begin{abstract}
It is known that unicellular LLT polynomials are related to the quasi-symmetric
chromatic polynomials of certain graphs by the $(t-1)$-transform of symmetric functions.
We investigate the extension of this transformation to various combinatorial
Hopf algebras and prove a noncommutative version of this property.
\end{abstract}

\maketitle

%%%%%%%%%%%%%%%%%%%%%%%%%%%%%%%%%%%%%%%%%%%%%%%%%%%%%%%%%%%%%%%%%%%%%%%%%%%%%%%
%%%%%%%%%%%%%%%%%%%%%%%%%%%%%%%%%%%%%%%%%%%%%%%%%%%%%%%%%%%%%%%%%%%%%%%%%%%%%%%
%%%%%%%%%%%%%%%%%%%%%%%%%%%%%%%%%%%%%%%%%%%%%%%%%%%%%%%%%%%%%%%%%%%%%%%%%%%%%%%
\section{Introduction}

In their proof of the shuffle conjecture \cite{CM}, Carlsson and Mellit obtain
a remarkable relation between unicellular LLT polynomials and the
quasi-symmetric chromatic polynomials \cite{SW} of certain graphs, namely
\begin{equation}\label{eq:X2LLT}
\X_G(t,X) = (t-1)^{-n}\,\LLT_G(t,(t-1)X).
\end{equation}
The graphs $G$ are simple undirected graphs with vertices labelled $1,\ldots,n$,
characterized by the property that if there is an edge $\{i,j\}$ with $i<j$,
then all the $\{i',j'\}$ with $i\leq i'<j'\leq j$ are also edges of $G$.
The number of such graphs is the Catalan number $c_n$. These are the
incomparability graphs of certain posets $P$, known as unit interval orders \cite{SW}.

Let $V(G)$ and $E(G)$ denote respectively the sets of vertices and edges of $G$.
A coloring of $G$ is a map $c:\ V(G)\rightarrow \NN^*$, which can be
identified with a word $c_1c_2\cdots c_n$.
A coloring is proper if $c_i\not=c_j$ whenever $\{i,j\}\in E(G)$. We denote by
$C(G)$ the set of proper colorings of $G$.
The chromatic quasi-symmetric function of $G$ expands in the $M$ basis of
$\QSym$ as~\cite{SW}
\begin{equation}
\X_G(t,X) = \sum_{c\in C(G)}t^{\asc_G(c)}x_{c_1}x_{c_2}\cdots x_{c_n}
        =  \sum_{c\in \PC(G)}t^{\asc_G(c)}M_{\ev(c)}(X),
\end{equation}
where $\PC(G)$ denotes the set of proper {\it packed colorings}, $\asc_G(c)$ is the
number of edges $\{i<j\}$ such that $c_i<c_j$, and $\ev(c)$ is the evaluation (or content)
of $c$, that is, the composition recording the number of occurences of each value of $c$.
It can be shown that for the above graphs, $\X_G(t)$ is actually a symmetric
function \cite{SW}.

On another hand, LLT polynomials are $t$-analogues of products of skew Schur
functions \cite{LLT1,HHL}.
By interpreting $s_1$ as $s_{\lambda/\mu}$ in various ways, one may obtain
different $t$-analogues of the characteristic $s_1^n$ of the regular
representation of $\SG_n$.
These can be parametrized by the same graphs as above, and their expression
given by the dinv statistic of Haglund, Haiman and Loehr~\cite{HHL} can be
rephrased as
\begin{equation}
\LLT_G(t,X) = \sum_{u\in PW_n}t^{\asc_G(u)}M_{\ev(u)}(X),
\end{equation}
where $u$ runs now over {\emph all} packed words of length $n$, regarded as
colorings of $G$.

Therefore, Equation \eqref{eq:X2LLT} tells us that the transformation $(t-1)X$
just eliminates the improper colorings, a fact which is far from obvious. A
``pedestrian'' proof can be found in the appendix of \cite{HX}.
The aim of this note is to provide a conceptual explanation (and a
generalization) in terms of combinatorial Hopf algebras.

%%%%%%%%%%%%%%%%%%%%%%%%%%%%%%%%%%%%%%%%%%%%%%%%%%%%%%%%%%%%%%%%%%%%%%%%%%%%%%%
%%%%%%%%%%%%%%%%%%%%%%%%%%%%%%%%%%%%%%%%%%%%%%%%%%%%%%%%%%%%%%%%%%%%%%%%%%%%%%%
%%%%%%%%%%%%%%%%%%%%%%%%%%%%%%%%%%%%%%%%%%%%%%%%%%%%%%%%%%%%%%%%%%%%%%%%%%%%%%%
\section{Word quasi-symmetric functions}

Let $A=\{a_1<a_2<\dots\}$ be a totally ordered alphabet.
The \emph{packed word} $u=\pack(w)$ associated with a word $w\in A^*$ is
a word over the alphabet of positive integers,
obtained by the following process. If $b_1<b_2<\ldots <b_r$ are the letters
occuring in $w$, $u$ is the image of $w$ by the homomorphism
$b_i\mapsto i$.

A word of positive integers $u$ is said to be \emph{packed} if $\pack(u)=u$. We denote by $\PW$ the
set of packed words.
With such a word, we associate the ``polynomial''
\begin{equation}
\M_u :=\sum_{\pack(w)=u}w\,.
\end{equation}
For example, with $A=\{1<2<3<4<5\}$, 
\begin{equation}
\begin{split}
\M_{13132} = &\ \ \ \ 13132 + 14142 + 14143 + 24243 \\
             &      + 15152 + 15153 + 25253 + 15154 + 25254 + 35354.
\end{split}
\end{equation}

The \emph{evaluation} $\ev(w)$ of a word $w$ is the sequence whose $i$-th
term $|w|_{a_i}$ is the number of times the letter $a_i$ occurs in $w$,
regarded as a finite integer vector by removing the trailing zeros.

Let $\KK$ be a field of characteristic $0$, assumed to contain all formal
series in the formal parameter $t$ used in the sequel.

Under the abelianization $\chi:\ \KK\langle A\rangle\rightarrow\KK[X]$, the
$\M_u$ are mapped to the monomial quasi-symmetric functions 
\begin{equation}
M_I := \sum_{j_1<j_2\ldots<j_r}x_{j_1}^{i_1}x_{j_2}^{i_2}\cdots x_{j_r}^{i_r},
\end{equation}
where $I=(|u|_a)_{a\in A}=(i_1,\ldots,i_r)$ is the evaluation vector of $u$.

The $\M_u$ span a subalgebra of $\KK\langle A\rangle$, called $\WQSym$
for Word Quasi-Symmetric functions, consisting of the invariants of the
noncommutative version of Hivert's quasi-symmetrizing action \cite{NCSF7},
which is defined by $\sigma\cdot w = w'$ where $w'$ is such that
$\std(w')=\std(w)$ and $\chi(w')=\sigma\cdot\chi(w)$, where $\std$ stands for
the usual standardization algorithm, namely the algorithm that sends any word
to the permutation having the same inversions.
Hence, two words are in the same $\SG(A)$-orbit iff they have the same packed
word.

When $A$ is infinite, $\KK\<A\>$ is interpreted as the algebra of formal
series of bounded degree. Exactly as in the case of symmetric or
quasi-symmetric functions, $\WQSym$ acquires then the structure of a Hopf
algebra, with the natural coproduct given by the ordinal sum of mutually
commuting alphabets.

The coproduct $A+B$ is indeed well-defined on $\WQSym$ and allows to consider
its graded dual $\WQSym^*$.
We shall denote by $\N_u\in\WQSym^*$ the dual basis of $\M_u$. 

The  algebra $\FQSym$ (Free Quasi-symmetric functions) may be defined as
the subalgebra of $\KK\<A\>$ spanned by the
\begin{equation}
\G_\sigma(A)= \sum_{\std(w)=\sigma}w,
\end{equation}
where $\sigma$ runs over all permutations. 
It is also a Hopf algebra for the same coproduct. It is self-dual, and the
dual basis $\F_\sigma=\G_\sigma^*$ can be identified with $\G_{\sigma^{-1}}$.
It is isomorphic to the Malvenuto-Reutenauer Hopf algebra \cite{MR,NCSF6}.

There is therefore an inclusion of Hopf algebras
$\iota:\ \FQSym\hookrightarrow\WQSym$ given by
\begin{equation}
\G_\sigma = \sum_{\std(u)=\sigma}\M_u,
\end{equation}
whose dual is the projection $\iota^*:\ \WQSym^*\twoheadrightarrow\FQSym$
\begin{equation}
\N_u\mapsto \F_{\std(u)}.
\end{equation}

Let $AB$ be the alphabet $\{ab|a\in A,\ b\in B\}$ endowed with the
lexicographic order on the pairs $(a,b)$.
It is easy to check that the coproduct $\delta:\ f\mapsto f(AB)$ is also
well-defined: packing with respect to the lexicographic order makes sense, and
\begin{equation}
\M_u(AB) = \sum_{\pack{\binom{v}{w}}=u}\M_v(A)\M_w(B),
\end{equation}
where $\binom{u}{v}$ denotes the word in \emph{biletters} $\binom{u_i}{v_i}$,
lexicographically ordered with priority to the top letter.
 
The dual of the coproduct $AB$ is an internal product
on each homogeneous component $\WQSym^*_n$ given by
\begin{equation}
\N_u * \N_v = \N_{\pack{\binom{u}{v}}},
\end{equation}
In this way,  $\WQSym_n^*$ gets identified with the Solomon-Tits
algebra of $\SG_n$ \cite{NT06} for all $n$.

There is also a Hopf embedding of $\Sym$ (Noncommutative symmetric functions)
into $\WQSym^*$ given by $S_n\mapsto \hat S_n := \N_{1^n}$:
\begin{equation}
S^I \mapsto \hat S^I=  \sum_{\ev(u)=I}\N_u.
\end{equation}
For example, $\hat S^{21}=\N_{112}+\N_{121}+\N_{211}$.

Under the projection $\iota^*$,
\begin{equation}
\hat S^I \mapsto \iota^*(\N_{1^{i_1}}\cdots\N_{1^{i_r}})
  = \F_{12\ldots i_1}\cdots\F_{12\ldots i_r}=S^I=\sum_{\Des(\sigma)\subseteq\Des(I)}\G_{\sigma},
\end{equation}
where $\Des(\sigma)$ denotes the descent set of $\sigma$, and $\Des(I)$ the set
encoded by the composition $I$.

This projection is compatible with the internal products: on $\Sym$, the
internal product is defined as dual to the coproduct $XY$ of $QSym$
\cite{Ge84,MR}.
By definition, it maps $f(AB)$ to $f(XY)$.
On $\FQSym$, the internal product on the $\F$-basis is ordinary composition:
$\F_\sigma*\F_\tau=\F_{\sigma\circ\tau}$, so that on the $\G$-basis,
$\G_\sigma*\G_\tau=\G_{\tau\circ\sigma}$.

Now, although the $*$ product of $\WQSym^*$ does not coincide with composition
on permutations, we have the following compatibility.

\begin{lemma}
\label{lem:perm}
Define a right action of $\SG_n$ on $\WQSym_n^*$ by
\begin{equation}
\N_u\cdot\sigma := \N_{u\sigma},\
     \text{where}\ u\sigma = u_{\sigma(1)}u_{\sigma(2)}\cdots u_{\sigma(n)}.
\end{equation}
Then, for any $I\vDash n$
\begin{equation}
\N_{u\sigma}*\hat S^I = (\N_u*\hat S^I)\cdot \sigma.
\end{equation}
\end{lemma}

\Proof If $\hat S^I$ contains $\N_v$, it contains  $\N_{v\tau}$ for all
permutations $\tau$, and
\begin{equation}\label{eq:permpack}
\pack{\binom{u\tau}{v}} = \pack{\binom{u}{v\tau^{-1}}}\cdot\tau.
\end{equation}
\qed

For example, with $u=111122$, $v=212211$, $\tau=451623$,
we have $u\tau = 121211$, $v\tau^{-1}=211212$, 
$\pack{\binom{121211}{212211}} = 232411$,
$\pack{\binom{111122}{211212}}=211234$,
and $211234\tau = 232411$.

This implies that $(f\cdot\sigma)*g =(f*g)\cdot\sigma$ for all
$f\in\WQSym^*_n, g\in\Sym_n$ and $\sigma \in\SG_n$.

\begin{remark}{\rm
A similar argument actually proves the existence of the descent algebra. 
If $\sigma=\std(v)$,
\begin{align}
\iota^*(\hat S^I*\N_v)
&=  \sum_{\ev(u)=I}\iota^*\left(\N_{\pack{\binom{u}{v}}}\right)\\
&=  \sum_{\ev(u)=I}\F_{\std{\binom{u}{v}}}\\
&=  \sum_{\ev(u)=I}\F_{\std{\binom{u\sigma^{-1}}{v\sigma^{-1}}}\circ\sigma}\\
&=  \sum_{\ev(u)=I}\F_{\std(u\sigma^{-1})\circ\sigma}\\
&=  S^I * \F_\sigma\\ 
&=  \iota^*(\hat S^I)*\iota^*(\N_v).
\end{align}
This implies in particular that in $\FQSym$,
\begin{equation}
S^I * S^J = \iota^*(\hat S^I*\hat S^J) = \sum_{M\in {\rm Mat}(I,J)} S^M,\\
\end{equation}
where ${\rm Mat}(I,J)$ denotes the set of nonnegative integer matrices with
row-sums vector $I$ and column-sums vector $J$ (cf. \cite{NCSF1}).
Hence, the $S^I=\sum_{\Des(\sigma^{-1})\subseteq\Des(I)}\F_{\sigma}$ span a sub $*$-algebra of
$\FQSym$ isomorphic to $(\Sym,*)$, which is therefore anti-isomorphic to the
Solomon descent algebra.
}
\end{remark}

%%%%%%%%%%%%%%%%%%%%%%%%%%%%%%%%%%%%%%%%%%%%%%%%%%%%%%%%%%%%%%%%%%%%%%%%%%%%%%%
%%%%%%%%%%%%%%%%%%%%%%%%%%%%%%%%%%%%%%%%%%%%%%%%%%%%%%%%%%%%%%%%%%%%%%%%%%%%%%%
%%%%%%%%%%%%%%%%%%%%%%%%%%%%%%%%%%%%%%%%%%%%%%%%%%%%%%%%%%%%%%%%%%%%%%%%%%%%%%%
\section{Transformations of alphabets}

%%%%%%%%%%%%%%%%%%%%%%%%%%%%%%%%%%%%%%%%%%%%%%%%%%%%%%%%%%%%%%%%%%%%%%%%%%%%%%%
\subsection{Transformations in $\QSym$}

Recall that the classical Cauchy identities for symmetric functions can be
extended to the dual pair of Hopf algebras $(QSym,\Sym)$ as follows. Let $X$
be a totally ordered alphabet of commutative variables, and $A$ be an alphabet
of noncommuting variables, also totally ordered.
The product alphabet $XA$ is the set of products $xa$ endowed with the
lexicographic order on the pairs $(x,a)$. We can thus define the
noncommutative symmetric functions of $XA$ and we have
\begin{equation}
\sigma_1(XA)
= \prod_{x\in X}^\rightarrow\prod_{a\in A}^\rightarrow(1-xa)^{-1}
= \sum_I M_I(X)S^I(A)
= \sum_I U_I(X)V_I(A)
\end{equation}
for any pair $(U,V)$ of mutually dual bases \cite{NCSF1} (the arrows mean that
the products are to be taken in increasing order).

We can introduce a second commutative alphabet $T$, and compute in two ways
\begin{equation}
\sigma_1(XTA) = \sum_I M_I(XT)S^I(A) = \sum_I M_I(X)S^I(TA).
\end{equation}

The alphabet $T$ denoted by $\frac1{1-t}$ is $\{t^n|n\ge 0\}$, ordered by
$t^i<t^j$ iff $i>j$ (which would be true for numerical values of $t$ such that
the geometric series converge). We introduce the notations
\begin{equation}
XT = \frac{X|}{|1-t}  \quad\text{and}\quad
TA = \frac{|A}{1-t|}.
\end{equation} 

The maps 
\begin{equation}
M_I(X)\mapsto M_I\left(\frac{X|}{|1-t}\right)
\quad\text{and}\quad
S^I(A) \mapsto S^I\left( \frac{|A}{1-t|}\right)
\end{equation} 
are algebra automorphisms, and their inverses are consistently denoted by
\begin{equation}
M_I(X)\mapsto M_I(X(1-t))  \quad\text{and}\quad
S^I(A)\mapsto S^I\left((1-t)A\right)
\end{equation} 
Here, $(1-t)$ is an example of a virtual alphabet. 
More generally, a virtual alphabet $T$ is defined as a morphism of algebras
$\chi_T$
\begin{equation}
M_I \mapsto  \chi_T(M_I) =: M_I(T)
\end{equation}
from $QSym$ to some commutative algebra. This defines the $TA$ transform in
$\Sym$ by
\begin{equation}
S_n(TA) := \chi_T(S_n(XA)) = \sum_{I\vDash n}M_I(T)S^I(A)
\end{equation}
and by duality, the $XT$ transform on $QSym$
\begin{equation}
M_I(XT) = \sum_J\<M_I(X),S^J(TA)\>M_J(X),
\end{equation}
where
\begin{align}
\<M_I(X),S^J(TA)\>&
  =\<\Delta^sM_I, S^{j_1}(TA)\otimes\cdots\otimes S^{j_s}(TA)\>\\
&=
\begin{cases}
M_{I_1}(T)\cdots M_{I_s}(T)&
 \text{if $I=I_1\cdots I_s$ with $I_k\vDash j_k$}\\
0 &
 \text{otherwise.}
\end{cases}
\end{align}

The $(t-1)$ transform is defined by writing $t-1=t(1-t^{-1})$ so that
$M_I(X(t-1))=t^{|I|}M_I(X(1-t^{-1}))$. 

We note for further reference the specializations 
\begin{align}
M_I\left(\frac{1}{1-t}\right)
&= \frac{t^{\maj(I)}}{(1-t^{i_1})(1-t^{i_1+i_2})\cdots(1-t^{i_1+\cdots+i_r})},
\ (\maj(I)=\sum_{d\in\Des(I)}d)\\
M_I\left(\frac{1}{t-1}\right)
&= \frac{1}{(1-t^{i_1})(1-t^{i_1+i_2})\cdots(1-t^{i_1+\cdots+i_r})},\\
M_I(1-t) & = (-1)^{\ell(I)-1}(t^{n-i_1}-t^n) \quad (I\vDash n),\\
M_I(t-1) & =  (-1)^{\ell(I)-1}(t^{i_1}-1).
\end{align}

%%%%%%%%%%%%%%%%%%%%%%%%%%%%%%%%%%%%%%%%%%%%%%%%%%%%%%%%%%%%%%%%%%%%%%%%%%%%%%%
\subsection{Transformations in $\WQSym$}

The $1/(1-t)$ transform may be extended from $\QSym$ to $\WQSym$ by setting
\begin{equation}
\frac{A|}{|1-t}= \{a_it^j|i\ge1,j\ge 0\}
\end{equation}
endowed with the total order $a_it^j<a_kt^l \Leftrightarrow i<k$ or $i=k$ and
$j>l$.
Then, the commutative image of $\M_u\left(\frac{A|}{|1-t}\right)$ is
$M_I\left(\frac{X|}{|1-t}\right)$, where $I=\ev(u)$.

The inverse transformations are consistently denoted by
$\M_u\mapsto \M_u(A(1-t))$ on $\WQSym$ and
$S^I\mapsto S^I((1-t)A)$ on $\Sym$.
These have been investigated in \cite{NTsuper,NCSF2}.

The adjoint map of  $\M_u\mapsto \M_u(A(1-t))$ is
$\N_u\mapsto \N_u*\sigma_1((1-t)A)$, and similarly for the inverse maps.

Although there is no known polynomial realization of $\WQSym^*$, it will be
convenient to define $\N_u(TA)$ as
\begin{align}
\N_u(TA) &:= \N_u*\sigma_1(TA) = \sum_v \M_v(T) \N_u*\N_v\\
         & = \sum_w\left(\sum_{v;\pack{\binom{u}{v}}=w}\M_v(T)\right)\N_w.
\end{align}

Let
\begin{equation}
V(u,w) = \{v|\pack {\binom{u}{v}}=w\}.
\end{equation}

\begin{proposition}
Let $u$ and $w$ be two packed words of the same size.
Let $w^{(i)}=\pack(w_{j_1}w_{j_2}\cdots w_{j_p})$,
where $\{j_1,\ldots,j_p\}=\{j|u_j=i\}$.
Then,
\begin{equation}
\sum_{v\in V(u,w)}\M_v = \M_{w^{(1)}}\M_{w^{(2)}}\cdots \M_{w^{(\max(u))}}.
\end{equation}
\end{proposition}

\Proof
Since the packing process commutes with the right action of the symmetric
group (see \eqref{eq:permpack}), we can apply to $u$ the smallest permutation
$\sigma$ such that $u\sigma$ is nondecreasing (i.e., $\sigma=\std(u)^{-1}$),
so that $\pack{\binom{u\sigma}{v\sigma}}=w\sigma$. We can therefore assume
that $u$ is nondecreasing.
First note that no relation is required between the letters of $v$
corresponding to different letters of $u$. The only order constraints are
among places where $u$ has identical letters, and these are the same as in the
corresponding letters of $w$. This is precisely the definition of the
convolution on packed words, describing the product of the $\M$
basis~\cite{NT06}.
\qed

Thus,
\begin{equation}
\N_u(TA) =
\sum_w ( \M_{w^{(1)}}\M_{w^{(2)}}\cdots \M_{w^{(\max(u))}})(T) \N_w,
\end{equation}
and by duality,
\begin{equation}
\label{eq:MuAT}
\M_u(AT) =
 \sum_v\left(\sum_{w\in V(v,u)}\M_w(T)\right)\M_v(A) =
 \sum_v(\M_{u^{(1)}}\cdots \M_{u^{(r)}})(T)\M_v(A).
\end{equation}

The morphism $\chi_T$ defining a virtual alphabet is naturally extended to
$\WQSym$ by setting $\M_u(T)=M_I(T)$, where $I=\ev(u)$.

A packed word $v$ is said to refine $u$ if for all $i<j$,
$v_i>v_j \Longleftrightarrow u_i\ge u_j$
and $v_i=v_j\Longrightarrow u_i=u_j$.
In this case, we write $v\raff u$.
This is the usual notion of refinement on set compositions: each block of $u$
is a union of \emph{consecutive} blocks of $v$.

For example, the packed words finer than $212$ are $212$, $213$, and $312$.
The packed words finer than $2122$ are
\begin{equation}
\label{raff2122}
2122, 2123, 2132, 3122, 2133, 3123, 3132, 2134, 2143, 3124, 3142, 4123, 4132.
\end{equation}

\begin{lemma}
\label{lem-wqsym2qsym}
The coefficient of $\M_v(A)$ in \eqref{eq:MuAT} is 0 if $u$ is not
finer than $v$, and equal to the coefficient of $M_{\ev(u)}(X)$ in
$M_{\ev(v)}(XT)$ otherwise.
\end{lemma}

\Proof
By definition, the words $u^{(i)}$ exist only when $u$ is finer than $v$,
and then, $ \M_{w^{(1)}}\M_{w^{(2)}}\cdots \M_{w^{(\max(u))}})(T)$
is equal to $M_{I_1}(T)\cdots M_{I_s}(T)$, where $I=\ev(v)$, $J=\ev(u)$
and $I_k\vDash j_k$ for all $k$.
\qed

%%%%%%%%%%%%%%%%%%%%%%%%%%%%%%%%%%%%%%%%%%%%%%%%%%%%%%%%%%%%%%%%%%%%%%%%%%%%%%%
%%%%%%%%%%%%%%%%%%%%%%%%%%%%%%%%%%%%%%%%%%%%%%%%%%%%%%%%%%%%%%%%%%%%%%%%%%%%%%%
%%%%%%%%%%%%%%%%%%%%%%%%%%%%%%%%%%%%%%%%%%%%%%%%%%%%%%%%%%%%%%%%%%%%%%%%%%%%%%%
\section{Dyck graphs}

\begin{definition}
A \em{Dyck graph} is a simple undirected graph $G$ with vertex set $V(G)=[n]$
and edge set $E(G)$ represented as pairs $(i<j)$ such that if $(i,j)\in E(G)$,
then $(i',j')\in E(G)$ for all $i\leq i'<j'\leq j$.
\end{definition}

Define for $\sigma\in\SG_n$
\begin{align}
%\inv_G(\sigma) & = \#\{(i<j)\in E(G)| \sigma^{-1}_i>\sigma^{-1}_j\}\\
\inv_G(\sigma) & = \#\{(i<j)\in E(G)| \text{$i$ is to the right of $j$
                                            in $\sigma$} \}\\
\Des_G(\sigma) & = \{i\,|\,\sigma_i>\sigma_{i+1}\ \text{and}\
                      (\sigma_{i+1},\sigma_i)\not\in E(G)\}\\
\maj_G(\sigma) &= \sum_{i\in\Des_G(\sigma)}i.
\end{align}

We shall make use of \emph{descent bottoms} of a permutation associated with
a graph, that are the values $\sigma_{i+1}$ such that
$\sigma_i>\sigma_{i+1}$ and $(\sigma_{i+1},\sigma_i)\not\in E(G)$.
%%%%%%%%%
For example, if $G$ is the graph
\begin{equation}
\begin{tikzpicture}
\begin{scope}[every node/.style={circle,scale=.5,fill=white,draw}]
    \node (A) at (0,0) {};
    \node (B) at (1*\taille,0) {};
    \node (C) at (2*\taille,0) {};
    \node (D) at (3*\taille,0) {};
    \node (E) at (4*\taille,0) {};
\end{scope}

\begin{scope}[>={Stealth[black]},
              every edge/.style={draw=black,thick}]
    \path [-] (A) edge (B);
    \path [-] (B) edge (C);
    \path [-] (C) edge (D);
    \path [-] (D) edge (E);
    \path [-] (C) edge[bend left=60] (E);
\end{scope}
\end{tikzpicture}
\end{equation}
(labelled 1--5 from left to right), and $\sigma=35142$, then
$\inv_G(\sigma)=\{(2,3),(4,5)\}$, $\Des_G(\sigma)=\{2,4\}$, $\maj_G(\sigma)=6$
and the descent bottoms of $\sigma$ are $\{1,2\}$.

\medskip
Set $\st_G(\sigma)=\inv_G(\sigma)+\maj_G(\sigma)$ \cite{SW,Kas}.
Recall the notation $[n]_t=1+t+\cdots t^{n-1}$.  
%%%%%%%
\begin{theorem}
Let $G$ be a Dyck graph. For any $\sigma\in\SG_{n-1}$, 
\begin{equation}
\sum_{\tau\in\sigma\shuffle n}t^{\st_G(\tau)}=[n]_t\, t^{\st_H(\sigma)},
\end{equation}
where $H$ is the restriction of $G$ to the interval $[1, n-1]$ and $\sigma\shuffle n$ means
the set of all words $\tau$ such that the restriction of $\tau$ to $[1,n-1]$
gives back $\sigma$.
\end{theorem}

\begin{corollary}
Let $\sigma|_{[1,k]}$ denote the restriction of $\sigma$ to the interval $[1,k]$.
The map
\begin{equation}\label{eq:stcod}
c(\sigma) = (\st_G(\sigma|_{[1,n-i]}))_{i=0\ldots n-1}
\end{equation}
is a bijection from $\SG_n$ to the set of integer vectors $v\in\NN^n$
such that $v_i\le n-i$.
\end{corollary}

Thus, $c$ is a code which interpolates between the Lehmer code (complete graph)
and the majcode (no edges).

In particular, we recover a particular case of a result of
Kasraoui~\cite{Kas}:
\begin{corollary}
For any Dyck graph $G$, the statistic $\st_G$ is Mahonian: 
\begin{equation}
\sum_{\sigma\in\SG_n}t^{\st_G(\sigma)}=[n]_t!
\end{equation}
\end{corollary}

The previous theorem is a direct consequence of the following lemma.

\begin{lemma}
\label{n:order}
Consider the permutations $\tau$ obtained from $\sigma$ by inserting $n$ at
each of the $n$ possible positions.

Then, $\st_G(\tau)-\st_H(\sigma)$ takes all the values from $0$ to $n-1$ in
this order if one visits the insertion positions in the following order: 

start with the rightmost position, then, from right to left, insert $n$ to the
left of the values $k$ such that $(k,n)\in E$ or that are descent bottoms of
$\sigma$, then from left to right, run through the remaining ones.
\end{lemma}

\Proof
First note that \eqref{eq:stcod} implies that $\st_G$ is Mahonian by induction
on $n$.

Let us now prove it.
Let $H$ be the restriction of $G$ to $[n-1]$. 
There are four cases to be distinguished according to whether $\tau$ is
obtained by inserting $n$:

\begin{enumerate}
\item at the end of $\sigma$: then $\st_G(\tau)=\st_H(\sigma)$.
\item to the left of a $k$ such that $(k,n)\in E(G)$. Then, $(k,k')\in E(G)$
for all $k<k'<n$, so that $k$ cannot be a descent bottom of $\sigma$.
Thus,
\begin{equation*}
\st_G(\tau)=\st_H(\sigma)+d_H(k)+e_G(k),
\end{equation*}
where $d_H(k)$ is the number of $H$-descent bottoms of $\sigma$  to the
right of $n$ (since each descent is shifted by one position to the right),
and $e_G(k)$ is the number of $k'$ to the right of $n$ such that
$(k',n)\in E(G)$ (since all these values have an inversion with $n$).
\item to the left of an $H$-descent bottom $k$. Then, the letter $\ell$
preceding $k$ in $\sigma$ is such that $(k,\ell)\not\in E(G)$, so that
$(l,k)\not\in E(G)$.
Therefore, inserting $n$ between $\ell$ and $k$ creates a descent $(n,k)$ in
$\tau$ which takes the place of the descent $(k,i)$ in $\sigma$ just one
position to the right.
Thus,
\begin{equation*}
\st_G(\tau)=\st_H(\sigma)+d_H(k)+e_G(k),
\end{equation*}
as in the previous case. Here, the $+1$ due to moving the descent bottom $k$
one place to the right is taken into account in $d_H(k)$.
\item to the left of a $k$ such that $(k,n)\not\in E(G)$ and $k$ is not an
$H$-descent bottom.
Then $(n,k)$ is a new descent and
\begin{equation*}
\st_G(\tau)=\st_H(\sigma)+\sigma^{-1}_k+d_H(k)+e_G(k).
\end{equation*}
\end{enumerate}

Let us now consider the sequence of insertions described in
Lemma~\ref{n:order}.

In the first part of the sequence going from right to left, one easily sees
that the values of $d_H(k)+e_G(k)$ increase by one at each step since we stop
at each element either creating an inversion with $n$ or being a descent
bottom.
In the second  part moving from left to right, the same property holds: at
each step, the value of $\sigma^{-1}_k+d_H(k)+e_G(k)$ increases by one since,
between two elements, $\sigma^{-1}_k$ changes by one plus the number of values
between these that are either descent bottoms or related with $n$ in $G$,
which is compensated by the fact that $d_H$ and $e_G$ decrease respectively on
descent bottoms or values connected to $n$ in $G$.

Finally, it is easily checked that both ways of evaluating the increment of
$\st_H$ corrresponding to the leftmost insertion position do agree, whence the
claim.
\qed 

\begin{note}{\rm
This argument is  similar to the one used for the maj-code 
(see, e.g., \cite{HNT06}). In particular, the definition of the sequence going
backward then forward to visit each possible insertion position of $n$ is
essentially the same: one just has to add a special case when $(k,n)\in E(G)$.
}
\end{note}

\begin{note}{\rm
The statistic $s_G(\sigma)$ interpolates in a Catalan number of ways between
the inversions number ($G$ is the complete graph) and the major index ($G$ is
the  graph with no edges).
}
\end{note}

\begin{example}{\rm 
Consider the graph
$\gsixex$
on $V(G)=[6]$ with $E(G)=\{12,23,24,34,45,46,56\}$
and the permutation $\sigma=52314$.

Then $4$ and $5$ belong to case (2), $1$ and $2$ belong to case (3), and $3$
belongs to case (4). Then the order is ${}^4 5^3 2^5 3^2 1^1 4^0$, where the
exponents encode the sequence.
}
\end{example}

\begin{proposition}\label{prop:dyck}
For a Dyck graph,
\begin{equation}
\X_G\left(t, \frac{1}{t-1}\right)=\frac{1}{(t-1)^n}.
\end{equation}
\end{proposition}
\Proof
According to \cite[Theorem 9.3]{SW}, the principal specialization of $\X_G$
satisfies
\begin{equation}
(q;q)_n\omega\X_G\left(t,\frac1{1-q}\right) =
 \sum_{\sigma\in\SG_n}t^{\inv_G(\sigma)}q^{\maj_G(\sigma)},
\end{equation}
where $(q;q)_n=(1-q)(1-q^2°\cdots (1-q^n)$.
For $q=t$, $\X_G$ being symmetric, this yields
\begin{equation}
(-1)^n(t;t)_n\X_G\left(t,\frac1{t-1}\right) =
 (t;t)_n\omega\X_G\left(t,\frac1{1-t}\right) = 
 \sum_{\sigma\in\SG_n}t^{\inv_G(\sigma)+\maj_G(\sigma)}=[n]_t!,
\end{equation}
where the first equality comes from the fact that for a symmetric function $f$ homogeneous
of degree $n$, $f(-X)=(-1)^n\omega f(X)$.
Dividing by $[n]_t!$, we are left with
\begin{equation}
(t-1)^n\,\X_G\left(t,\frac1{t-1}\right) =  1.\qed
\end{equation}
%\qed

\goodbreak

%%%%%%%%%%%%%%%%%%%%%%%%%%%%%%%%%%%%%%%%%%%%%%%%%%%%%%%%%%%%%%%%%%%%%%%%%%%%%%%
%%%%%%%%%%%%%%%%%%%%%%%%%%%%%%%%%%%%%%%%%%%%%%%%%%%%%%%%%%%%%%%%%%%%%%%%%%%%%%%
%%%%%%%%%%%%%%%%%%%%%%%%%%%%%%%%%%%%%%%%%%%%%%%%%%%%%%%%%%%%%%%%%%%%%%%%%%%%%%%
\section{The Guay-Paquet Hopf algebra}

In his proof of the Shareshian-Wachs conjecture \cite{GP}, Guay-Paquet
introduces a Hopf algebra $\GPG$ based on ordered graphs, depending on a
parameter $t$, and such that the map sending a graph to itsromatic quasi-symmetric
function is a morphism of Hopf algebras $\GPG\rightarrow QSym$.

Its basis consists of finite simple undirected graphs with vertices labelled
by the integers from 1 to $n=|V(G)|$. The product is the shifted concatenation:
$G\cdot H = G\cup H[n]$ where $H[n]$ is $H$ with labels shifted by the number
$n$ of vertices of $G$.

The parameter $t$ arises in the coproduct. If $G$ is a graph on $n$ vertices
and $w\in[r]^n$, regarded as a coloring of $G$, we denote by $G|_w$ the tensor
product $G_1\otimes\cdots\otimes G_r$ of the restrictions of $G$ to vertices
colored $1,2,\ldots,r$. The $r$-fold coproduct is then
\begin{equation}\label{eq:coprod}
\Delta^r G = \sum_{w\in[r]^n} t^{\asc_G(w)} G|_w.
\end{equation}

At $t=1$, $\GPG$ becomes cocommutative and is isomorphic to an algebra
introduced in~\cite{Schm}.

It is also proved in \cite{GP} that the subspace $\GPD$ of $\GPG$ spanned by
Dyck graphs is a Hopf subalgebra. At $t=1$, it is a free cocommutative graded
connected Hopf algebra of graded dimension Catalan, and is therefore
isomorphic to $\CQSym$~\cite{HNT05}.

From these properties, we obtain a simple conceptual proof
of~\eqref{eq:X2LLT}:

\begin{proposition}
\begin{equation}
(t-1)^n \X_G\left(t, \frac{X|}{|t-1} \right)
  =  \sum_{u\in PW_n}t^{\asc_G(u)}M_{\ev(u)}(X)
  = \LLT_G(t,X).
\end{equation}
\end{proposition}

\Proof
If $G$ is a Dyck graph and $u$ a packed word, denote by $G_i(u)$ the
restriction of $G$ to the vertices $j$ such that $u_j=i$.
Then the coefficient of $M_I(X)$ in $\X_G\left(t, \frac{X|}{|t-1} \right)$ is
\begin{equation}
\left\<\X_G\left(t, \frac{X|}{|t-1} \right),S^I\right\>
=\left\<\X_G(t,X),S^I\left(\frac{|A}{t-1|}\right) \right\>.
\end{equation}
Dualizing the product $S^I=S_{i_1}\cdots S_{i_r}$, this is equal to
\begin{equation}
\left\<\Delta^r\X_G(t,X),
         (S_{i_1}\otimes\cdots\otimes S_{i_r})\left(\frac{|A}{t-1|}\right)
  \right\>,
\end{equation}
and since $G\mapsto \X_G$ is a morphism of Hopf algebras, the iterated coproduct
can be evaluated by \eqref{eq:coprod}, which yields
\begin{equation}
\sum_{u\in\PW_n} t^{\asc_G(u)}
                    \prod_i \X_{G_i(u)}\left(t,\frac{1}{t-1}\right)
	= (t-1)^{-n}\sum_{u\in\PW_n}t^{\asc_G(u)},
\end{equation}
by Prop. \ref{prop:dyck}. \qed

Thanks to Lemma \ref{lem-wqsym2qsym}, this argument can be extended to the
noncommutative case.

%ICI
%%%%%%%%%%%%%%%%%%%%%%%%%%%%%%%%%%%%%%%%%%%%%%%%%%%%%%%%%%%%%%%%%%%%%%%%%%%%%%%
%%%%%%%%%%%%%%%%%%%%%%%%%%%%%%%%%%%%%%%%%%%%%%%%%%%%%%%%%%%%%%%%%%%%%%%%%%%%%%%
%%%%%%%%%%%%%%%%%%%%%%%%%%%%%%%%%%%%%%%%%%%%%%%%%%%%%%%%%%%%%%%%%%%%%%%%%%%%%%%
\section{The noncommutative chromatic quasi-symmetric function}

\subsection{A noncommutative analogue of $\X_G$}

Given a Dyck graph $G$, define
\begin{equation}
\XX_G(t,A) =  \sum_{c\in\PC(G)}t^{\asc_G(c)}\M_{c}(A) \in \WQSym.
\end{equation}

For example,
\begin{equation}
\XX_{\left(\gunun\right)} =  \M_{1}.
\end{equation}
\begin{equation}
\XX_{\left(\gdeuxun\right)} = \M_{11} + \M_{12} + \M_{21},
\end{equation}
\begin{equation}
\XX_{\left(\gdeuxdeux\right)} = t\,\M_{12} + \M_{21}.
\end{equation}
\begin{equation}
\XX_{\left(\gtroisun\right)} = \sum_{w\in PW_3} \M_{w},
\end{equation}
\begin{equation}
\begin{split}
\XX_{\left(\gtroisdeux\right)}
 =& t\,\M_{121} + t\,\M_{122} + t\,\M_{123} + t\,\M_{132}
  + \M_{211} \\
 &+ \M_{212} + \M_{213} + t\,\M_{231}
  + \M_{312} + \M_{321},
\end{split}
\end{equation}
\begin{equation}
\begin{split}
\XX_{\left(\gtroistrois\right)}
 =& t\,\M_{112} + \M_{121} + t\,\M_{123} + \M_{132}
 + t\,\M_{212} \\
 &+ t\,\M_{213} + \M_{221} + \M_{231}
 + t\,\M_{312} + \M_{321},
\end{split}
\end{equation}
\begin{equation}
\XX_{\left(\gtroisquatre\right)}
 = t\,\M_{121} + t^2\,\M_{123} + t\,\M_{132} + t\,\M_{212}
 + t\,\M_{213} + t\,\M_{231} + t\,\M_{312} + \M_{321},
\end{equation}

\begin{equation}
\XX_{\left(\gtroiscinq\right)}
 = t^3\,\M_{123} + t^2\,\M_{132} + t^2\,\M_{213} + t\,\M_{231}
 + t\,\M_{312} + \M_{321}.
\end{equation}

\begin{proposition}
$G\mapsto \XX_G(A)$ is a morphism of Hopf algebras from $\GPG$ to $\WQSym$.
\end{proposition}

\Proof
The argument is essentially the same as for $QSym$.
Multiplicativity is clear:
\begin{align}
\XX_{G_1}\XX_{G_2}
& = \sum_{(u_1,u_2)\in\PC(G_1)\times\PC(G_2)}
       t^{\asc_{G_1}(u_1)+\asc_{G_2}(u_2)}\M_{u_1}\M_{u_2}\\
& = \sum_{(u_1,u_2)\in\PC(G_1)\times\PC(G_2)}
       t^{\asc_{G_1}(u_1)+\asc_{G_2}(u_2)}
       \sum_{\ogf{v=v_1v_2}{\ogf{\pack(v_1)=u_1}{\pack(v_2)=u_2}}}\M_v\\
& = \sum_{v\in\PC(G_1G_2)}t^{\asc_{G_1G_2}(v)}\M_v.
\end{align}
Next, the coefficient of $\M_{u_1}\otimes\M_{u_2}$ in
$\Delta \XX_G$ is nonzero if and only if $u_1\in\PC(G_1)$ and
$u_2\in\PC(G_2)$ for some splitting of the vertices of $G$ into two
complementary subsets, detemined by a word $w\in\{1,2\}^n$ as in \eqref{eq:coprod}.

Each such splitting determines a proper coloring $u$ of $G$: color the
vertices of $G_1$ with $u_1$ and those of $G_2$ with the shifted word
$u_2[\max(u_1)]$ (recall that the shifted word $u[k]$ is obtained by adding
$k$ to all values of $u$). Thus, $u_1$ and $u_2$ are the restrictions of $u$
to two consecutive intervals, so that $\M_{u_1}\otimes \M_{u_2}$ occurs in
$\Delta \M_u$. In particular, note that $u$ belongs to the shifted shuffle of
$u_1$ and $u_2$.

Conversely, any $u\in\PC(G)$ and any term $\M_{u_1}\otimes\M_{u_2}$ occuring
in $\Delta\M_u$ uniquely determines a splitting of $V(G)$ into two
complementary subsets: the vertices of $G_1$ correspond to the positions of
the letters of the subword $u_1$ of $u$. This proves that at $t=1$,
$\XX$ is indeed a morphism of coalgebras.

Now, in the above situation, we have
\begin{equation}
\asc_G(u)=\asc_{G_1}(u_1)+\asc_{G_2}(u_2)+r
\end{equation}
where $r$ is the number of  edges $(i<j)$ of $G$ with $u_i<u_j$  which are neither in $G_1$
nor in $G_2$. These correspond precisely to the $G$-ascents of the word
$w\in\{1,2\}^n$ determining the splitting.
\qed

From the product rule of the $\M$ basis of $\WQSym$, one can easily
check that
\begin{equation}
\XX_{\left(\gunun\right)} \XX_{\left(\gdeuxun\right)} 
= \XX_{\left(\gtroisun\right)},
\end{equation}
and that
\begin{equation}
\XX_{\left(\gunun\right)} 
\XX_{\left(\gdeuxdeux\right)}
= \XX_{\left(\gtroistrois\right)}.
\end{equation}

For the coproduct, one can check the following example:
\begin{equation}
\begin{split}
\Delta\XX_{\left(\gtroisdeux\right)}
&= \XX_{\left(\gtroisdeux\right)} \otimes 1 +
\XX_{\left(\gunun\right)} \otimes
   \left((1+t)\XX_{\left(\gdeuxun\right)} + \XX_{\left(\gdeuxdeux\right)}
   \right) \\
&+ \left((1+t)\XX_{\left(\gdeuxun\right)} + \XX_{\left(\gdeuxdeux\right)}
   \right)
\otimes \XX_{\left(\gunun\right)} +
1\otimes \XX_{\left(\gtroisdeux\right)}.
\end{split}
\end{equation}

This example also allows to check that the restriction of the coproduct to the
subalgebra of Dyck graphs is cocommutative.

\subsection{Noncommutative unicellular LLT polynomials}

\begin{theorem}
\begin{equation}
(t-1)^n \XX_G\left(t, \frac{A|}{|t-1} \right)
 = \sum_{u\in PW_n}t^{\asc_G(u)}\M_{u}(A).
\end{equation}
The r.h.s. is therefore a noncommutative lift of the LLT polynomial $\LLT_G$.
\end{theorem}

\Proof
The coefficient of $\M_v$ in $\XX_G\left(t,\frac{A|}{|t-1}\right)$ is
\begin{equation}\label{eq:MMinXX}
c_v(t) =
 \left\<S^{\ev(v)},\sum_{\ogf{u\in\PC(G)}{u\raff v}}t^{\asc_G(u)}
                   M_{\ev(u)} \left(\frac{X|}{|t-1}\right)\right\>.
\end{equation}
Up to a power of $t$, the sum in the right-hand side of the bracket is the
product of the quasi-symmetric chromatic polynomials of the graphs $G_i(v)$
evaluated at $\frac{X|}{|t-1}$. The power of $t$ corresponds to the
$G$-ascents of $v$ on the deleted edges, that is
\begin{equation}
c_v(t) =
 t^{\asc_G(v)}\prod_i X_{G_i(v)}\left(\frac{1}{t-1}\right)
= \frac{t^{\asc_G(v)}}{(t-1)^n}.
\end{equation}
\qed

\begin{definition}
Given a Dyck graph $G$,
the \emph{non-commutative LLT polynomial} $\bLLT_G$ is
\begin{equation}
\bLLT_G := \sum_{u\in PW_n}t^{\asc_G(u)}\M_{u}(A).
\end{equation}
\end{definition}
\begin{note}{\rm
Alternatively, the r.h.s of \eqref{eq:MMinXX} can be interpreted as a duality
bracket for the pair $(\WQSym^*,\WQSym)$:
\begin{align}
c_v(t) &=
 \left\<S^{\ev(v)},\sum_{\ogf{u\in\PC(G)}{u\raff v}}t^{\asc_G(u)}
                   \M_{u} \left(\frac{A|}{|t-1}\right)\right\>\\
&= \left\<\N_{i_1}\cdots \N_{i_r},
          \sum_{\ogf{u\in\PC(G)}{u\raff v}}t^{\asc_G(u)}
                   \M_{u} \left(\frac{A|}{|t-1}\right)\right\>\\
&= \left\<\N_{i_1}\otimes\cdots \otimes\N_{i_r},
          \sum_{\ogf{u\in\PC(G)}{u\raff v}}t^{\asc_G(u)}
                   \Delta^r\M_{u} \left(\frac{A|}{|t-1}\right)\right\>
\end{align}
and evaluating the iterated coproducts $\Delta^r\M_{u}$ leads to the same
conclusion.
}
\end{note}

%For example, with $G=\gcinqun$ (chain graph) and $v=22112$, the minimal proper
%colorings of $G$ finer than $v$ are
%\begin{equation}
%34124, 34214, 43124, 43214, 34123, 34213, 43123, 43213.
%\end{equation}
%These are the colorings of $G_1(v)$ with its minimal number (two) of colors
%$1,2$ and of $G_2(v)$ with its miminal number (two again) of colors, hence
%colors $3,4$.
%
%Then the coefficient of $22112$ in $\bLLT_G$ 

%%%%%%%%%%%%%%%%%%%%%%%%%%%%%%%%%%%%%%%%%%%%%%%%%%%%%%%%%%%%%%%%%%%%%%%%%%%%%%%
\subsection{Special case: path graphs}

Let $G_n$ be the graph on $[n]$ with edges $(i,i+1)$. Then,
\begin{equation}
\bLLT_{G_n} = \sum_{u\in\PW_n}t^{\asc_{G_n}(u)}\M_u.
\end{equation}
If we embed $\Sym$ in $\WQSym$ by sending $S_n$ to the sum of nonincreasing
words
\begin{equation}
S_n\mapsto
   \sum_{u\in\PW_n,\ u\downarrow} \M_u\quad
\Leftrightarrow \ \Lambda_n\mapsto \M_{12\cdots n},
\end{equation}
we can write
\begin{equation}
\bLLT_{G_n} = \sum_{w\in A^n}t^{\asc(w)} w
           = \sum_{I\vDash n}(t-1)^{n-\ell(I)}\Lambda^I
\end{equation}
so that
\begin{equation}
\XX_{G_n} = \sum_{I\vDash n}\frac{\Lambda^I(A(t-1))}{(t-1)^{\ell(I)}}
\end{equation}
which gives back by commutative image the generating series of \cite[C.2]{SW}.
The images of $S_n$ and $\Lambda_n$ by the $A\mapsto A(t-1)$ transform are
given by
\begin{equation}
\sigma_1(A(t-1))
= \sum_{u\downarrow}\M_u(A(t-1))
= \sum_{u\downarrow} t^{|u|-\max(u)}(t-1)^{\max(u)}\M_u
= \prod_{i\ge 1}^\leftarrow \frac{1-a_i}{1-ta_i}
\end{equation}
and its inverse ($\lambda_{-t}:=(\sigma_t)^{-1}$)
\begin{equation}
\lambda_{-1}(A(t-1))
= \sum_{u=12\cdots n}\M_u(A(t-1))
= \sum_{u\uparrow} (1-t)^{\max(u)}\M_u
= \prod_{i\ge 1}^\rightarrow \frac{1-ta_i}{1-a_i}.
\end{equation}
Hence,
\begin{equation}
\left(\sum_{n\ge 0}\XX_{G_n}\right)^{-1}
= 1+\sum_{n\ge 1}(-1)^n\sum_{u\uparrow}(1-t)^{\max(u)-1}\M_u.
\end{equation}
At $t=1$, this gives back the well-known fact that the sum of all Smirnov
words is the inverse of the alternating sum of constant words.

%%%%%%%%%%%%%%%%%%%%%%%%%%%%%%%%%%%%%%%%%%%%%%%%%%%%%%%%%%%%%%%%%%%%%%%%%%%%%%%
%%%%%%%%%%%%%%%%%%%%%%%%%%%%%%%%%%%%%%%%%%%%%%%%%%%%%%%%%%%%%%%%%%%%%%%%%%%%%%%
%%%%%%%%%%%%%%%%%%%%%%%%%%%%%%%%%%%%%%%%%%%%%%%%%%%%%%%%%%%%%%%%%%%%%%%%%%%%%%%
\section{The Dyck graphs subalgebra of $\WQSym$}

The goal of this section is to prove

\begin{theorem}
The restriction of the morphism of Hopf algebras $G\mapsto \XX_G(t,A)$ from
$\GPG$ to $\WQSym$ to the subalgebra $\GPD$ of Dyck graphs is injective.
\end{theorem}

We shall prove that the images of the Dyck graphs are already linearly
independent for $t=1$.

%%%%%%%%%%%%%%%%%%%%%%%%%%%%%%%%%%%%%%%%%%%%%%%%%%%%%%%%%%%%%%%%%%%%%%%%%%%%%%%
\subsection{The Hopf algebra $\WSym$}

The $\XX_G(1,A)$ are the noncommutative chromatic polynomials defined by
Gebhard \cite{Geb,GebSa}, and thus belong to the algebra of symmetric
functions in noncommuting variables $a_i$, denoted here by $\WSym$.

It consists of the invariants of $\SG(A)$ acting by automorphisms
on the free algebra $\K\<A\>$.
Two words $u=u_1\cdots u_n$ and $v=v_1\cdots v_n$ are in the same orbit
whenever $u_i=u_j \Leftrightarrow v_i=v_j$. Thus, orbits are parametrized by
set partitions into at most $|A|$ blocks. Assuming that $A$ is infinite, we
obtain an algebra based on all set partitions whose monomial basis is defined
by
\begin{equation}
\m_\pi(A)=\sum_{w\in O_\pi} w
\end{equation}
where $O_\pi$ is the set of words such that $w_i=w_j$ iff $i$ and $j$ are
in the same block of $\pi$.

As an example of expansion of a chromatic polynomial in terms of the $\m$, we
have
\begin{equation}
\XX_{(\gtroisdeux)}(1) = \m_{121} + \m_{122} + \m_{123}.
\end{equation}

The product of the $\m$ is given by the rule
\begin{equation}
\m_{\pi'} \m_{\pi''} = \sum_{\pi\in E(\pi',\pi'')} \m_\pi,
\end{equation}
where $E(\pi',\pi'')$ consists of all set partitions whose parts are either a
part of $\pi'$, a part of $\pi''$, or a union of a part of $\pi'$ and a part
of $\pi''$.

Since set partitions are equivalence classes of set compositions which are in
bijection with packed words, we shall often denote a set partition as the
lexicographically minimal packed word in its class, which amounts
to representing a set partition by the set composition obtained by ordering
its blocks w.r.t. their minima. For example, $\{\{1,4\}\{2,5\},\{3\}\}$
will be represented by $12312$.

Let us illustrate this notation on two examples of the product:
\begin{equation}
\m_{1}\m_{1123} = \m_{12234} + \m_{11123} + \m_{12213} + \m_{12231},
\end{equation}
and
\begin{equation}
\m_{1123}\m_1 = \m_{11234} + \m_{11233} + \m_{11232} + \m_{11231}.
\end{equation}

%%%%%%%%%%%%%%%%%%%%%%%%%%%%%%%%%%%%%%%%%%%%%%%%%%%%%%%%%%%%%%%%%%%%%%%%%%%%%%%
\subsection{The chromatic polynomials}

To prove the linear independence of the images of the Dyck graphs, we shall
show that they are triangular with respect to a basis of a subalgebra of
$\WSym$ based on \emph{nonnesting partitions}.

Define the \emph{denesting} $\dn(\pi)$ of a set partition $\pi$ as the
nonnesting partition $\pi'$ obtained by iterating the following process: for
each sequence $i<j<k<l$ such that $j$ and $k$ are in a block $B_1$ of $\pi$
and $i$ and $l$ in another block $B_2$ containing no intermediate value $i<r<l$,
i.e., $B_2=\{m_1<\cdots<m_p=i<m_{p+1}=l<\cdots m_r\}$,
split $B_2$ into $B'_2=\{m_1,\ldots,i\}$ and $B''_2=\{l,\ldots<m_r\}$.

Up to $n=3$, all set partitions are fixed by the denesting algorithm and there
is only one set partition $\pi$ of size $4$ such that $\dn(\pi)\not=\pi$. In
terms of set partitions, it is $\{\{1,4\},\{2,3\}\}$ and
$\dn(\pi)=\{\{1\},\{2,3\},\{4\}\}$.
In terms of packed words, it is $1221$ and $\dn(1221)=1223$.

If $\pi=12341312$, then $\dn(\pi)=12341356$.

\begin{proposition}
For a nonnesting partition $\pi$, define
\begin{equation} 
\mt_\pi = \sum_{\dn(\pi')=\pi}\m_{\pi'}.
\end{equation}
Then, the $\mt_\pi$ form the basis of a subalgebra of $\WSym$ of homogeneous
dimensions given by the Catalan numbers.
\end{proposition}

For example,
\begin{equation} 
\mt_{1223} = \m_{1223} + \m_{1221}.
\end{equation}
\begin{equation}
\mt_{12334} = \m_{12334} + \m_{12331} + \m_{12332}.
\end{equation}
\begin{equation}
\mt_{12233} = \m_{12233} + \m_{12211}.
\end{equation}
\begin{equation}
\mt_{12324} = \m_{12324} + \mt_{12321}.
\end{equation}

\Proof
Since the product of the $\m$-basis is multiplicity-free, we just have to
check that for any partition $\pi'$, $\m_{\pi'}$ occurs in
$\mt_{\pi_1}\mt_{\pi_2}$ if and only if $\m_{\dn(\pi')}$ occurs in this
product.

If $\pi_1\vdash [n_1]$ and $\pi_2\vdash [n_2]$, then  
$\m_{\pi'}$ occurs in $\mt_{\pi_1}\mt_{\pi_2}$ if and only
if $\dn(\pi'|_{[1,n_1]})=\pi_1$ and  $\dn(\pi'|_{[n_1,n_1+n_2]})=\pi_2$.

Since the denesting process is obviously compatible with restriction to
intervals, this is equivalent to  $\dn(\pi')|_{[1,n_1]}=\pi_1$ and
$\dn(\pi')|_{[n_1,n_1+n_2]}=\pi_2$, which is the condition for
$\m_{\dn(\pi')}$ to occur in $\mt_{\pi_1}\mt_{\pi_2}$.
\qed

\begin{lemma}
For a Dyck graph $G$,
\begin{equation}
\XX_G(1,A)
 = \sum_{\pi(i)\not=\pi(j) \text{\ if\ } (i,j)\in E(G)} \mt_\pi
\end{equation}
where the sum runs over nonnesting partitions, and $\pi(i)$ denotes the block containing $i$.
\end{lemma}

\Proof
If $\m_{\pi'}$ occurs in $\XX_G$, then so does $\mt_{\dn(\pi')}$, since
$\dn(\pi')$ is finer than $\pi'$, so is associated with proper colorings as
well.

Conversely, if $\m_{\pi'}$ does not occur in $\XX_G$, there exist $i,j$ with
$|j-i|$ minimal such that $(i,j)\in E(G)$ and $i,j$ in the same block of
$\pi'$. By minimality of $|j-i|$, $i$ and $j$ are consecutive in their block.
Moreover, $(i,j)\in E(G)$ implies that $(i',j')\in E(G)$ for all $i<i'<j'<j$.
Still by minimality of $|j-i|$, $i+1,\ldots,j-1$ are all in different blocks.
Hence, $i$ and $j$ would not be separated by the denesting process, so that
$\m_{\dn(\pi')}$ does not occur in $\XX_G$ either.
\qed

There is a simple bijection $\eta$ between nonnesting partitions $\pi$
(represented as diagrams of arcs) and Young diagrams $\lambda$ contained in
the staircase partition $(n-1,\ldots,2,1)$,
represented as sets of cells above the diagonal in an $n\times n$ square: 
the arcs of $\pi$ are the
coordinates of the  corners of $\lambda$.
For example, the partition $\lambda=(221)$ corresponds to the nonnesting
partition $i\{1,3\},\{2,4\},\{5\}\}$ which is read on the coordinates of the corners of the diagram.
The edges of corresponding graph $G$ are the coordinates of the empty cells above the diagonal,
$(1,2),(2,3),(3,4),(3,5),(4,5)$. 
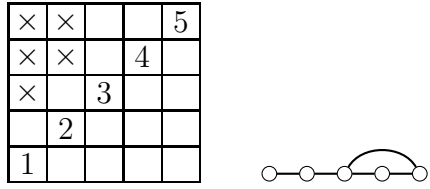
\begin{equation*}
\young{\times &\times & & & 5\cr
       \times & \times & &4& \cr
       \times & & 3 &&\cr
        & 2 &&&\cr
        1&&&&\cr}
\quad\quad
\begin{tikzpicture}
\begin{scope}[every node/.style={circle,scale=.5,fill=white,draw}]
    \node (A) at (0,0) {};
    \node (B) at (1*\taille,0) {};
    \node (C) at (2*\taille,0) {};
    \node (D) at (3*\taille,0) {};
    \node (E) at (4*\taille,0) {};
\end{scope}

\begin{scope}[>={Stealth[black]},
              every edge/.style={draw=black,thick}]
    \path [-] (A) edge (B);
    \path [-] (B) edge (C);
    \path [-] (C) edge (D);
    \path [-] (D) edge (E);
    \path [-] (C) edge[bend left=60] (E);
\end{scope}
\end{tikzpicture}
\end{equation*}

Thanks to that bijection, there is a natural partial order on nonnesting
partitions: the Young lattice restricted to partitions contained in the
staircase. We shall say that $\pi'\le\pi$ if the image of $\pi'$ is included
in the image of $\pi$.

\begin{lemma}
Given a Dyck graph $G$,
\begin{equation}
\XX_G(1,A) = \sum_{\pi'\le \pi_G}\mt_{\pi'},
\end{equation}
where the sum runs over nonnesting partitions smaller than the nonnesting
partition $\pi_G$ corresponding to the Young diagram encoding $G$.
\end{lemma}

\Proof
Let $\lambda=\eta(\pi_G)$.
Let $\pi'$ be a nonnesting set partition. If $\eta(\pi')\subseteq\lambda$
then, thanks to the bijection between partitions and Dyck graphs, for all
$(i,j)\in E(G)$, all pairs $(i',j')$ such that $i\leq i'<j'\leq j$ are also
edges of $G$, so that $i'$ and $j'$ can never be in the same part of $\pi$,
hence of $\pi'$. So thanks to the previous lemma, $\mt_{\pi'}$ appears in the
expansion of $\XX_G$.

Conversely, if $\eta(\pi')\not\subseteq\lambda$, then there exists a corner
$(i,j)$ of $\eta(\pi')$ that does not belong to $\lambda$.
Then $(i,j)$ is an edge of $G$ since it is an empty cell in $\lambda$ but $i$
and $j$ are in the same block of $\pi'$ since they are consecutive by
definition. So $\mt_{\pi'}$ does not appear in the expansion of $\XX_G$.
\qed

For example,
\begin{equation}
\XX_{\left(\gunun\right)}(1) =  \mt_{1}.
\end{equation}
\begin{equation}
\XX_{\left(\gdeuxun\right)}(1) = \mt_{11} + \mt_{12},
\end{equation}
\begin{equation}
\XX_{\left(\gdeuxdeux\right)}(1) = \mt_{12}.
\end{equation}
\begin{equation}
\XX_{\left(\gtroisun\right)}(1) = \mt_{111} + \mt_{112} + \mt_{122} +
  \mt_{121} + \mt_{123},
\end{equation}
\begin{equation}
\XX_{\left(\gtroisdeux\right)}(1)
 = \mt_{122} + \mt_{121} + \mt_{123},
\end{equation}
\begin{equation}
\XX_{\left(\gtroistrois\right)}(1)
 = \mt_{112} + \mt_{121} + \mt_{123}
\end{equation}
\begin{equation}
\XX_{\left(\gtroisquatre\right)}(1)
 = \mt_{121} + \mt_{123}
\end{equation}
\begin{equation}
\XX_{\left(\gtroiscinq\right)}(1)
 = \mt_{123}
\end{equation}

%%%%%%%%%%%%%%%%%%%%%%%%%%%%%%%%%%%%%%%%%%%%%%%%%%%%%%%%%%%%%%%%%%%%%%%%%%%%%%%
%%%%%%%%%%%%%%%%%%%%%%%%%%%%%%%%%%%%%%%%%%%%%%%%%%%%%%%%%%%%%%%%%%%%%%%%%%%%%%%
%%%%%%%%%%%%%%%%%%%%%%%%%%%%%%%%%%%%%%%%%%%%%%%%%%%%%%%%%%%%%%%%%%%%%%%%%%%%%%%
\section{A multiplicative basis}

Recall that the \emph{reverse refinement order}, denoted by $\le$, on
compositions is such that $I=(i_1,\ldots,i_k)\ge J=(j_1,\ldots,j_l)$ iff
$\{i_1,i_1+i_2,\ldots,i_1+\cdots+i_k\}$ contains
$\{j_1,j_1+j_2,\ldots,j_1+\cdots+j_l\}$.
In this case, we say that $I$ is finer than $J$.
For example, $(2,1,2,3,1,2)\ge (3,2,6)$.

This order can be extended to packed words as follows.
To avoid confusion with the order of packed words defined previously,
we say that $w$ is \emph{strongly finer} than $w'$, and write $w\ge w'$, iff
$w$ and $w'$ have same standardized word and the evaluation of $w$ is finer
than the evaluation of $w'$. It also amounts to asking that $w$ and $w'$ have
same standardized word and that $w\raff w'$.

For example, the packed words strongly finer than $212$ are $212$ and $213$.
The packed words strongly finer than $2122$ are
\begin{equation}
2122, 2123, 2133, 2134.
\end{equation}
to be compared with the packed words finer that $2122$ shown
in~\eqref{raff2122}.

Given a permutation, let us define the set $A(\sigma)$ of the \emph{advances}
of $\sigma$ as the set of values $i$ such that $i+1$ is to the right of $i$ in
$\sigma$. Note that it is the complementary set over $[1,n-1]$ of the usual
recoils. Let $\DST(\sigma)$ be the set of packed words of standardization
$\sigma$.

\begin{lemma}
\label{destand}
Let $\sigma$ be a permutation. Then the elements of $\DST(\sigma)$
are in bijection with the subsets of $A(\sigma)$.

\noindent
In particular, this set of words has a natural structure of boolean lattice.
\end{lemma}

\Proof
Let $i$ be an element of $A(\sigma)$ and let $j<k$ be the respective
positions of $i$ and $i+1$ in $\sigma$.
Then in any element $w$ of $\DST(\sigma)$, either $u_j=u_k$ or
$u_k=u_j+1$. Since there are two independent choices for all elements of
$A(\sigma)$, the result holds.

The inverse map going from 
$\DST(\sigma)$ to subsets of $A(\sigma)$ is given by the rule:
put $i$ in its corresponding subset of $A(\sigma)$ if $u_k=u_j+1$.
\qed

\begin{corollary}
Let $w$ be a word. Then the elements strongly coarser than $w$ are an interval
of the boolean order of $\DST(w)$ described in Lemma~\ref{destand}, hence
themselves a boolean order consisting of the subsets of $A(\std(w))$ that,
following the notations of the previous lemma, necessarily contains the
elements $i$ such that $w_j=w_k$.
\end{corollary}

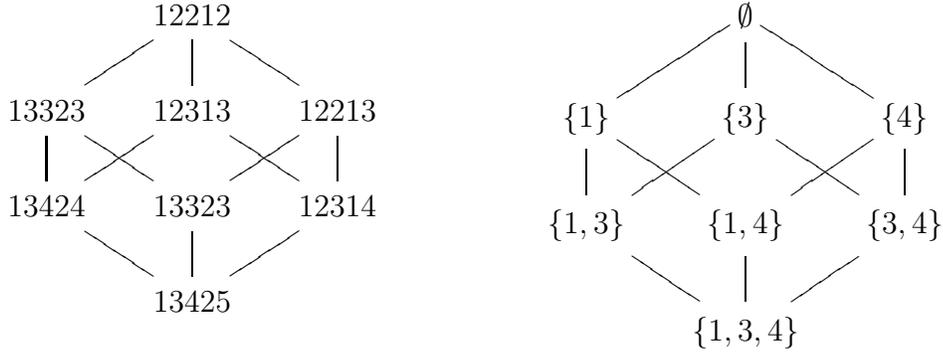
\begin{figure}[ht]
\newdimen\vcadre\vcadre=0.2cm % marges verticales de la boite
\newdimen\hcadre\hcadre=0.1cm % marges horizontales de la boite
$\xymatrix@R=0.6cm@C=7mm{
% niv 0
&& *{\GrTeXBox{12212}}\arx1[ld]\arx1[d]\arx1[rd]\\
% niv 1
&  *{\GrTeXBox{13323}}\arx1[d]\arx1[rd]
&  *{\GrTeXBox{12313}}\arx1[ld]\arx1[rd]
&  *{\GrTeXBox{12213}}\arx1[ld]\arx1[d]\\
% niv 2
&  *{\GrTeXBox{13424}}\arx1[rd]
&  *{\GrTeXBox{13323}}\arx1[d]
&  *{\GrTeXBox{12314}}\arx1[ld] \\
% niv 3
&&*{\GrTeXBox{13425}}
}$
\hskip1cm
$\xymatrix@R=0.6cm@C=7mm{
% niv 0
&& *{\GrTeXBox{\emptyset}}\arx1[ld]\arx1[d]\arx1[rd]\\
% niv 1
&  *{\GrTeXBox{ \{1\} }}\arx1[d]\arx1[rd]
&  *{\GrTeXBox{ \{3\} }}\arx1[ld]\arx1[rd]
&  *{\GrTeXBox{ \{4\} }}\arx1[ld]\arx1[d]\\
% niv 2
&  *{\GrTeXBox{ \{1,3\} }}\arx1[rd]
&  *{\GrTeXBox{ \{1,4\} }}\arx1[d]
&  *{\GrTeXBox{ \{3,4\} }}\arx1[ld] \\
% niv 3
&&*{\GrTeXBox{ \{1,3,4\} }}
}$
\caption{\label{fig-bool13425}%
The boolean orders of the packed words of standardized $13425$ and
the corresponding subsets of $\{1,3,4\}$.}
\end{figure}

The boolean order of the packed words of standardized $13425$ is given in
Figure~\ref{fig-bool13425}.

Let~\cite{NT06,BZ}
\begin{equation} 
\Phi_u := \sum_{v\ge u} \M_v.
\end{equation} 
For example,
\begin{equation}
\Phi_{111} = \M_{111} + \M_{112} + \M_{122} + \M_{123} ; \quad
%\Phi_{112} = \M_{112} + \M_{123} ;
\Phi_{212} = \M_{212} + \M_{213}.
\end{equation}
\begin{equation}
\Phi_{133142} = \M_{133142} + \M_{134152} + \M_{144253} + \M_{145263}.
\end{equation}

Since $(\Phi_u)$ is triangular over  $(\M_u)$, it is a basis of $\WQSym$. Note
that the order used for summation is a restriction of the refinement order on
compositions, so is a boolean lattice. 
Hence,
\begin{equation}
\M_u = \sum_{v\ge u} (-1)^{\max(v)-\max(u)} \Phi_v.
\end{equation}
For example,
\begin{equation}
\M_{133142} = \Phi_{133142} - \Phi_{134152} - \Phi_{144253} + \Phi_{145263}.
\end{equation}

By construction, the basis $\Phi$ satisfies a product formula similar to
that of Gessel's basis $F_I$ of $\QSym$ (whence the choice of notation).
We shall not state it since we will not need it in the sequel but here follows
an example illustrating the similarity with Gessel's basis. 

\begin{equation}
\begin{split}
\Phi_{1} \Phi_{121} &= \Phi_{1121}+\Phi_{2121}+\Phi_{3121}+\Phi_{2132},\\
F_1 F_{21} &= F_{31} + F_{22} + F_{211} + F_{121}.
\end{split}
\end{equation}

\begin{proposition}
\label{ncchrom2phi}
The noncommutative $t$-chromatic function is $\Phi$-positive:
\begin{equation}
\XX_G(t) = \sum_{\sigma\in\SG_n}t^{\asc_G(\sigma)}\Phi_{\min_G(\sigma)},
\end{equation}
where $\min_G(\sigma)$ is the packed word $u$ defined as follows: let 
\begin{equation}
S = \{i| \sigma^{-1}_{i-1}<\sigma^{-1}_i \text{\ and\ }
         (\sigma^{-1}_{i-1},\sigma^{-1}_i)\not\in E(G)\}.
\end{equation}
Then, 
\begin{equation}
u_i = \sigma_i - |S\cap [1,i]|.
\end{equation}
\end{proposition}

All non trivial examples (excluding the case of the complete graph where $S$
is always empty) of size $3$ are given below.

First, here are all sets $S$ and then all packed words $\min_G(\sigma)$.
\medskip
\begin{center}
\begin{tabular}{|c||c|c|c|c|c|c|}
\hline
  &123  &132  &213  &231  &312  &321 \\
\hline
\hline
$\gtroisun$    & $\{2,3 \}$ & $\{2 \}$ & $\{3 \}$ & $\{3 \}$ & $\{2 \}$ &
$\emptyset$\\
\hline
$\gtroisdeux$   & $\{3 \}$ & $\{2 \}$ & $\{3 \}$ & $\emptyset$ & $\{2 \}$ &
$\emptyset$\\
\hline
$\gtroistrois$  & $\{2 \}$ & $\{ 2\}$ & $\{3 \}$ & $\{3 \}$ & $\emptyset$ &
$\emptyset$\\
\hline
$\gtroisquatre$ & $\emptyset$ & $\{ 2\}$ & $\{ 3\}$ & $\emptyset$ &
$\emptyset$ & $\emptyset$\\
\hline
\end{tabular}
\end{center}

\medskip
\begin{center}
\begin{tabular}{|c||c|c|c|c|c|c|}
\hline
  &123  &132  &213  &231  &312  &321 \\
\hline
\hline
$\gtroisun$  & 111    & 121    & 212    & 221    & 211    & 321   \\
\hline
$\gtroisdeux$  & 122    & 121    & 212    &  231   & 211    & 321   \\
\hline
$\gtroistrois$  &  112   & 121    &   212  & 221    &  312   & 321   \\
\hline
$\gtroisquatre$  &   123  & 121    &  212   & 231    & 312    &  321  \\
\hline
\end{tabular}
\end{center}

We then deduce
\begin{equation}
\XX_{\left(\gunun\right)} =  \Phi_{1}.
\end{equation}
\begin{equation}
\XX_{\left(\gdeuxun\right)} = \Phi_{11} + \Phi_{21},
\end{equation}
\begin{equation}
\XX_{\left(\gdeuxdeux\right)} = t\,\Phi_{12} + \Phi_{21}.
\end{equation}

\begin{equation}
\XX_{\left(\gtroisun\right)}
 = \Phi_{111} + \Phi_{121} + \Phi_{212} + \Phi_{221} + \Phi_{211} +
   \Phi_{321}.
\end{equation}

\begin{equation}
\XX_{\left(\gtroisdeux\right)}
 = t\,\Phi_{122} + t\,\Phi_{121} + \Phi_{212} + t\,\Phi_{231}
  + \Phi_{211} + \Phi_{321},
\end{equation}

\begin{equation}
\XX_{\left(\gtroistrois\right)}
 = t\,\Phi_{112} + \Phi_{121} + t\,\Phi_{212} + \Phi_{221}
 + t\,\Phi_{312} + \Phi_{321},
\end{equation}
\begin{equation}
\XX_{\left(\gtroisquatre\right)}
 = t^2\,\Phi_{123} + t\,\Phi_{121} + t\,\Phi_{212} + t\,\Phi_{231}
 + t\,\Phi_{312} + \Phi_{321},
\end{equation}

\begin{equation}
\XX_{\left(\gtroiscinq\right)}
 = t^3\,\Phi_{123} + t^2\,\Phi_{132} + t^2\,\Phi_{213} + t\,\Phi_{231}
 + t\,\Phi_{312} + \Phi_{321}.
\end{equation}

\Proof
Since $\Phi_w$ expanded in the $\M$ basis is a sum over elements of
$\DST(\std(w))$, we have to show two things:
first, the monomials in $t$ that are coefficients of the $\M$ are constant
among elements of $\DST(\std(w))$ and among the elements of $\DST(\std(w))$,
the elements appearing in $\XX_G$ expanded in the $\M$ basis are exactly the
elements finer than $\min_G(\sigma)$.

\medskip
Concerning the coefficients (monomials in $t$) of these elements, if $w$
appears with a coefficient $t^i$ then this is also the coefficient of
$\std(w)$ and more generally of any element strongly finer than $w$. Indeed,
the power of $t$ counts the ascents among the pairs of edges of $G$ and
that does not change from $w$ to $w'\geq w$ if $w$ appears in $\XX_G$.
This comes from the fact that all ascents of $w$ are ascents of $w'$ and that
the only positions $(i,j)$ that could add an ascent from $w$ to $w'$ are
those such that $w'_i>w'_j$ and $w_i=w_j$, but in that case $(i,j)$ cannot
be an edge of $G$ since $w$ is a proper coloring of $G$, and hence cannot count
as an ascent of $w'$.

\medskip
Let us now show that the packed words with standardized $\sigma$
appearing in $\XX_G$ are the elements finer than $\min_G(\sigma)$.
Recall that thanks to Lemma~\ref{destand}, all words having a given
standardized $\sigma$ form a boolean lattice when equipped with the strong
refinement order, and that any element corresponds to a subset of the set of
values $i$ such that $i-1$ is to the left of $i$ in $\sigma$ (or, equivalently
$\sigma^{-1}_{i-1} < \sigma^{-1}_i$).

Since we are looking for the packed words that are proper colorings of $G$, it
is pretty clear that the subsets containing an $i$ such that
$(\sigma^{-1}_{i-1},\sigma^{-1}_i)\in E(G)$ cannot bring proper colorings
since two connected vertices would have the same color. Conversely, excluding
those values necessarily brings a proper coloring.

So all packed words appearing in the expansion of $\XX_G$ in the $M$ basis
with a given standardized word $\sigma$ are strongly finer than
$\min_G(\sigma)$ and have all same coefficient, whence the statement.
\qed

\begin{corollary}[\cite{SW}, Thm 3.1]
The chromatic quasi-symmetric function of a Dyck graph is $F$-positive
and its expansion is
\begin{equation}
X_{G}(t) = \sum_{\sigma\in\SG_n} t^{\inv_G(\sigma)}
              F_{\widetilde{DES_P(\sigma)}},
\end{equation}
where $\inv_G(\sigma)$ is the pairs $(i<j)$ such that
$(\sigma_i,\sigma_j)\in E(G)$ and $\sigma_i>\sigma_j$,
and $DES_P(\sigma)$ is the composition encoding the set of $i$ such that
$(\sigma_i,\sigma_{i+1})\not\in E(G)$ and $\sigma_i>\sigma_{i+1}$,
and where $\tilde{}$ denotes the conjugate composition.
\end{corollary}

\Proof
Our formula in $\WQSym$ is projected to this expression by the morphism
sending $\Phi_w$ to $F_{\ev(w)}$: the contribution of $\sigma$ in our
equation is the contribution of $\sigma'=\inv(r(\sigma))$ in their equation,
where $r(\sigma)$ sends each value $i$ to $n+1-i$ if $\sigma\in\SG_n$.

Indeed, our ascents of a permutation go to the inversions of~\cite{SW} through
$\inv\circ r$ since $r$ changes ascents to inversions and
$\inv$ exchanges values and positions. Moreover, two values $i$ and $i+1$ of
$\sigma$ are equal in $\min_G(\sigma)$ (hence are not a descent of the
composition $\ev(\min_G(\sigma))$) iff they are increasing and their positions
do not correspond to an edge of $G$. Since $\sigma'$ can also be described as
$\sigma'={\overline{\inv(\sigma)}}$, where $\bar w$ denotes the mirror-image of
$w$, this exactly translates in $\sigma'$ as the values in positions
$n-i$, $n+1-i$ that decrease and do not form an edge of $G$, which is exactly
the definition of the conjugate of $DES_P(\sigma')$.
\qed

For example, the graph $G=\gsixex$ and the permutation $314652$ contribute in
our case to a term $t^3\M_{212321}$, whereas $G$ and
$(\overline{314652})^{-1}=(453162)^{-1}=453162$ contribute in the case
of~\cite{SW} as $t^3 M_{231}$.

%%%%%%%%%%%%%%%%%%%%%%%%%%%%%%%%%%%%%%%%%%%%%%%%%%%%%%%%%%%%%%%%%%%%%%%%%%%%%%%
\subsection{Noncommutative LLT polynomials}

To formulate the noncommutative analogue of the $F$-expansion of unicellular
LLT polynomials, we need another lift of the $F$-basis, defined as
\begin{equation}
\vPhi_u = \sum_{v\ge \bar u} \M_{\bar v}
\end{equation}
where the bar involution sends a word to its mirror image.
Thus,
\begin{equation}
\vPhi_u  =    \overline{(\Phi_{\bar u})}
\end{equation}
and its commutative image is again
$\vPhi_u(X)=F_{\ev(\bar u)}(X)=F_{\ev(u)}(X)$.

For a permutation $\sigma$ and a Dyck graph $G$, define
\begin{equation}
\text{min}_G'(\sigma) := \overline{\text{min}_{\bar G}(\bar \sigma)},
\end{equation}
where $\bar G$ is the mirror image of $G$ (which amounts to relabeling
the vertices by $i\mapsto n+1-i$).

Here are the non-trivial examples of $\min'_G$ for $n=3$:

\medskip
\begin{center}
\begin{tabular}{|c|c|c|c|c|c|c|}
\hline
  &123  &132  &213  &231  &312  &321 \\
\hline
\hline
$\gtroisun$     & 123  & 122 & 112  & 121  & 212 & 111   \\
\hline
$\gtroisdeux$   & 123  & 122 & 213  & 121  & 212 & 211   \\
\hline
$\gtroistrois$  & 123  & 132 & 112  & 121  & 212 & 221   \\
\hline
$\gtroisquatre$ & 123  & 132 & 213  & 121  & 212 & 321  \\
\hline
\end{tabular}
\end{center}
\medskip

With these definitions, Proposition \ref{ncchrom2phi} translates into:
\begin{proposition}
The noncommutative $t$-chromatic polynomial is $\vPhi$-positive:
\begin{equation}
\XX_G(t) = \sum_{\sigma\in\SG_n}t^{\asc_G(\sigma)}\vPhi_{\min'_G(\sigma)}.
\end{equation}
\end{proposition}

For example,
\begin{equation}
\XX_{\left(\gtroisun\right)}
 = \vPhi_{123} + \vPhi_{122} + \vPhi_{112} + \vPhi_{121} + \vPhi_{212}
 + \vPhi_{111}.
\end{equation}
\begin{equation}
\XX_{\left(\gtroisdeux\right)}
 = t\,\vPhi_{123} + t\,\vPhi_{122} + \vPhi_{213} + t\,\vPhi_{121}
 + \vPhi_{212} + \vPhi_{211},
\end{equation}
\begin{equation}
\XX_{\left(\gtroistrois\right)}
 = t\,\vPhi_{123} + \vPhi_{132} + t\,\Phi_{112} + \vPhi_{121}
 + t\,\vPhi_{212} + \vPhi_{221},
\end{equation}
\begin{equation}
\XX_{\left(\gtroisquatre\right)}
 = t^2\,\vPhi_{123} + t\,\vPhi_{132} + t\,\vPhi_{212} + t\,\vPhi_{121}
   + t\,\vPhi_{212} + \vPhi_{321}.
\end{equation}

\begin{theorem}
The noncommutative unicellular LLT polynomials are $\vPhi$-positive:
\begin{equation}
\bLLT_G
 = \sum_{\sigma\in\SG_n}t^{\asc_G(\sigma)}\vPhi_{\min'_{G_\emptyset}(\sigma)}
\end{equation}
where $G_\emptyset$ is the graph with ($n$ vertices, $n$ omitted) and
no edges.
\end{theorem}

\Proof
If $\M_v$ occurs in $\vPhi_u$, then $\asc_G(v)=\asc_G(u)$ for any graph $G$.
Also, since $\min'_{G_\emptyset}(\sigma)$ is $\bar v$ where $v$ is the minimal
element of $\DST(\bar\sigma)$,
\begin{equation}
\sum_{\sigma\in\SG_n}\vPhi_{\min'_{G_\emptyset}(\sigma)}=\sum_{u\in\PW_n}\M_u,
\end{equation}
so that 
\begin{equation}
\bLLT_G = 
\sum_{u\in\PW_n}t^{\asc_G(u)}\M_u=
\sum_{\sigma\in\SG_n}t^{\asc_G(\sigma)}\vPhi_{\min'_{G_\emptyset}(\sigma)}.
\end{equation}
\qed

For example,
\begin{equation}
\bLLT_{\left(\gtroisun\right)}
 = \vPhi_{123} + \vPhi_{122} + \vPhi_{112} + \vPhi_{121}
 + \vPhi_{212} + \vPhi_{111},
\end{equation}
\begin{equation}
\bLLT_{\left(\gtroisdeux\right)}
 = t\,\vPhi_{123} + t\,\vPhi_{122} + \vPhi_{112} + t\,\vPhi_{121}
 + \vPhi_{212} + \vPhi_{111},
\end{equation}
\begin{equation}
\bLLT_{\left(\gtroistrois\right)}
 = t\,\vPhi_{123} + \vPhi_{122} + t\,\Phi_{112} + \vPhi_{121}
 + t\,\vPhi_{212} + \vPhi_{111},
\end{equation}
\begin{equation}
\bLLT_{\left(\gtroisquatre\right)}
 = t^2\,\vPhi_{123} + t\,\vPhi_{122} + t\,\vPhi_{112} + t\,\vPhi_{121}
 + t\,\vPhi_{212} + \vPhi_{111}.
\end{equation}

%%%%%%%%%%%%%%%%%%%%%%%%%%%%%%%%%%%%%%%%%%%%%%%%%%%%%%%%%%%%%%%%%%%%%%%%%%%%%%%
%%%%%%%%%%%%%%%%%%%%%%%%%%%%%%%%%%%%%%%%%%%%%%%%%%%%%%%%%%%%%%%%%%%%%%%%%%%%%%%
%%%%%%%%%%%%%%%%%%%%%%%%%%%%%%%%%%%%%%%%%%%%%%%%%%%%%%%%%%%%%%%%%%%%%%%%%%%%%%%

\end{document}